\documentclass[a4paper,12pt]{article}
\usepackage[dvips]{epsfig,graphicx}
\usepackage{amsmath,amsfonts,verbatim}

\newcommand{\RR}{{\mathbb R}}

\newcommand{\CC}{{\mathbb C}}

\newcommand{\qed}{\hfill\hbox{\rule{3pt}{6pt}}}
\newcommand{\proof}{{\sc Proof. }}

\newcommand{\G}{\Gamma}
\newcommand{\mat}{\hbox{Mat}}

\newcommand{\Es}{E^*}
\newcommand{\MX}{\mat_X(\CC)}

\newcommand{\la}{\langle}
\newcommand{\h}{\hat}
\newcommand{\ra}{\rangle}

\newtheorem{theorem}{Theorem}[section]
\newtheorem{lemma}[theorem]{Lemma}
\newtheorem{corollary}[theorem]{Corollary}

\newtheorem{definition}[theorem]{Definition}

\newtheorem{notation}[theorem]{Notation}

\title{ON BIPARTITE DISTANCE-REGULAR GRAPHS WITH EXACTLY TWO IRREDUCIBLE T-MODULES WITH ENDPOINT TWO}

\author{Mark S. MacLean \\
        Mathematics Department\\
        Seattle University \\
        901 Twelfth Avenue \\
        Seattle WA 98122-1090, USA \\
        macleanm@seattleu.edu \and
        \v{S}tefko Miklavi\v{c}\footnote{Supported in part by ARRS - Javna
agencija za raziskovalno dejavnost Republike Slovenija, program no. P1-0285.} \\ 
        University of Primorska \\
        Andrej Maru\v{s}i\v{c} Institute\\
        Muzejski trg 2 \\ 
        6000 Koper, Slovenia \\
        stefko.miklavic@upr.si}

\begin{document}
\maketitle

\begin{abstract}
  Let $\G$ denote a bipartite distance-regular graph with diameter $D \ge 4$ and valency $k \ge 3$. 
  Let $X$ denote the vertex set of $\G$, and let $A$ denote the adjacency matrix of $\G$. 
  For $x \in X$ let $T=T(x)$ denote the subalgebra of $\MX$ generated by $A, \Es_0, \Es_1, \ldots, \Es_D$, where for $0 \le i \le D$, $\Es_i$ represents the projection onto the $i$th subconstituent of $\G$ with respect to $x$.  We refer to $T$ as the {\em Terwilliger algebra} of $\G$ with respect to $x$.  An irreducible $T$-module $W$ is said to be {\em thin} whenever dim $\Es_i W \le 1$ for $0 \le i \le D$.  By the {\em endpoint} of $W$ we mean min$\{i | \Es_iW \ne 0\}$.  For $0 \le i \le D$, let $\G_i(z)$ denote the set of vertices in $X$ that are distance $i$ from vertex $z$.  
  Define a parameter $\Delta_2$ in terms of the intersection numbers by $\Delta_2 = (k-2)(c_3-1)-(c_2-1)p^2_{22}$.  In this paper we prove the following are equivalent:  
 (i) $\Delta_2>0$ and for $2 \le i \le D - 2$ there exist complex scalars 
$\alpha_i, \beta_i$ with the following property: for all $x, y, z \in X$ such that 
$\partial(x, y) = 2, \: \partial(x, z) = i, \: \partial(y, z) = i$ we have
$  \alpha_i + \beta_i  |\G_1(x) \cap \G_1(y) \cap \G_{i-1}(z)| = |\G_{i-1}(x) \cap \G_{i-1}(y) \cap \G_1(z)|;$
(ii) For all $x \in X$ there exist up to isomorphism
exactly two irreducible modules for the Terwilliger algebra $T(x)$
with endpoint two, and these modules are thin.
\end{abstract}

\section{Introduction}
\label{sec:intro}

\bigskip

In this paper we obtain a combinatorial characterization of bipartite distance-regular graphs with exactly two irreducible modules of the Terwilliger algebra of endpoint $2$, both of which are thin (see Section 2 for formal definitions). Our combinatorial characterization is closely related to the 2-homogeneous property of Curtin \cite{Cu0} and Nomura \cite{No1}.


Throughout this introduction let $\G$ denote a bipartite distance-regular graph with diameter $D \geq 
4$ and valency $k \ge 3$.  Let $X$ denote the vertex set of $\G$.  For $x 
\in X$, let $T=T(x)$ denote the Terwilliger algebra of $\G$ with respect to $x$.  It is known that there exists a unique irreducible $T$-module with endpoint 0, and this module is thin \cite[Proposition 8.4]{egge1}.  Moreover, Curtin showed that up to isomorphism $\G$ has exactly one irreducible $T$-module with endpoint 1, and this module is thin \cite[Corollary 7.7]{Cu1}.

We now discuss the irreducible $T$-modules of endpoint 2.
For $0 \le i \le D$, let $\G_i(z)$ denote the set of vertices in $X$ that are distance $i$ from vertex $z$.  
In \cite[Theorem 3.11]{Cu4}, Curtin proved that the following are equivalent:  
 (i) For all $i$ $(2 \le i \le D - 2)$ and for all $x, y, z \in X$ with $\partial(x,y)=2, \partial(x,z)=i, \partial(y,z)=i$, the number  $|\G_1(x) \cap \G_1(y) \cap \G_{i-1}(z)|$ is independent of $x,y,z;$
(ii) For all $x \in X$ there exists a unique irreducible $T$-module for the Terwilliger algebra $T(x)$ with endpoint 2, and this module is thin.
When these equivalent conditions hold, $\G$ is said to be {\em almost} 2-{\em homogeneous.}

  Now define a parameter $\Delta_2$ in terms of the intersection numbers by $\Delta_2 = (k-2)(c_3-1)-(c_2-1)p^2_{22}$.  In this paper we prove the following are equivalent:  
 (i) $\Delta_2>0$ and for $2 \le i \le D - 2$ there exist complex scalars 
$\alpha_i, \beta_i$ with the following property: for all $x, y, z \in X$ such that 
$\partial(x, y) = 2, \: \partial(x, z) = i, \: \partial(y, z) = i$ we have
$  \alpha_i + \beta_i  |\G_1(x) \cap \G_1(y) \cap \G_{i-1}(z)| = |\G_{i-1}(x) \cap \G_{i-1}(y) \cap \G_1(z)|;$
(ii) For all $x \in X$ there exist up to isomorphism
exactly two irreducible modules for the Terwilliger algebra $T(x)$
with endpoint two, and these modules are thin.
We also compute $\alpha_i, \beta_i$ in terms of the intersection numbers of $\Gamma$.

We remark that this paper is part of a continuing effort to understand and classify the bipartite distance-regular graphs with at most two irreducible modules of the Terwilliger algebra with endpoint 2, both of which are thin.  Please see \cite{Cu2}--\cite{Cu4}, \cite{MT}--\cite{Maclean2} for more work from this ongoing project.

\section{Preliminaries}
\label{sec:prelim}

In this section we review some definitions and basic results concerning distance-regular graphs. 
See the book of Brouwer, Cohen and Neumaier \cite{BCN} for more background information.

\medskip \noindent
Let $\CC$ denote the complex number field and let $X$ denote a nonempty finite set. Let $\MX$ denote 
the $\CC$-algebra consisting of all matrices whose rows and columns are indexed by $X$ and whose 
entries are in $\CC$. Let $V=\CC^X$ denote the vector space over $\CC$ consisting of column vectors 
whose coordinates are indexed by $X$ and whose entries are in $\CC$. We observe $\MX$ acts on $V$ by 
left multiplication. We call $V$ the {\it standard module}. We endow $V$ with the Hermitean inner 
product $\langle \, , \, \rangle$ that satisfies $\langle u,v \rangle = u^t\overline{v}$ for 
$u,v \in V$, where $t$ denotes transpose and $\overline{\phantom{v}}$ denotes complex conjugation. For 
$y \in X$ let $\hat{y}$ denote the element of $V$ with a 1 in the $y$ coordinate and 0 in all 
other coordinates. We observe $\{\hat{y}\;|\;y \in X\}$ is an orthonormal basis for $V$.

\medskip \noindent
Let $\G = (X,R)$ denote a finite, undirected, connected graph, without loops or multiple edges, 
with vertex set $X$ and edge set $R$. Let $\partial$ denote the path-length distance function for 
$\G$, and set $D := \mbox{max}\{\partial(x,y) \;|\; x,y \in X\}$. We call $D$  the {\it diameter} of 
$\G$. For a vertex $x \in X$ and an integer $i$ let $\G_i(x)$ denote the set of vertices at 
distance $i$ from $x$. We abbreviate $\G(x) = \G_1(x)$. For an integer $k \ge 0$ we say $\Gamma$ is 
{\em regular with valency} $k$ whenever $|\G(x)|=k$ for all $x \in X$. We say $\G$ is 
{\it distance-regular} whenever for all integers $h,i,j\;(0 \le h,i,j \le D)$ and for all vertices 
$x,y \in X$ with $\partial(x,y)=h,$ the number
\begin{eqnarray*}
p_{ij}^h = | \G_i(x) \cap \G_j(y)|
\end{eqnarray*}
is independent of $x$ and $y$. The $p_{ij}^h$ are called the {\it intersection numbers} of $\G$. 

\medskip \noindent
For the rest of this paper we assume $\Gamma$ is distance-regular with diameter $D \ge 4$. 
Note that $p^h_{ij}=p^h_{ji}$ for $0 \le h,i,j \le D$. For convenience set 
$c_i:=p_{1, i-1}^i \, (1 \le i \le D)$, $a_i:=p_{1i}^i \, (0 \le i \le D)$, 
$b_i:=p_{1, i+1}^i \, (0 \le i \le D-1)$, $k_i:=p_{ii}^0 \, (0 \le i \le D)$, and $c_0=b_D=0$. By the 
triangle inequality the following hold for $0 \le h,i,j \le D$: (i) $p_{ij}^h=0$ if one of $h,i,j$ is 
greater than the sum of the other two; (ii) $p_{ij}^h \ne 0$ if one of $h,i,j$ equals the sum of the 
other two. In particular $c_i \ne 0$ for $1 \le i \le D$ and $b_i \ne 0$ for $0 \le i \le D-1$. We 
observe that $\G$ is regular with valency $k=k_1=b_0$ and that 
\begin{equation}
\label{kk}
c_i+a_i+b_i=k \qquad \qquad (0 \le i \le D).
\end{equation}
Note that $k_i = |\G_i(x)|$ for $x \in X$ and $0 \le i \le D$. By \cite[p.~127]{BCN},
\begin{equation}
\label{k}
  k_i = \frac{b_0 b_1 \cdots b_{i-1}}{c_1 c_2 \cdots c_i} \qquad \qquad (0 \le i \le D).
\end{equation}

\medskip \noindent 
We recall the Bose-Mesner algebra of $\G$. For $0 \le i \le D$ let $A_i$ denote the matrix in 
$\MX$ with $(x,y)$-entry
\begin{equation}
\label{dm1}
  (A_i)_{x y} = \left\{ \begin{array}{ll}
                 1 & \hbox{if } \; \partial(x,y)=i, \\
                 0 & \hbox{if } \; \partial(x,y) \ne i \end{array} \right. \qquad (x,y \in X).
\end{equation}
For notational convenience, we define $A_{D+1}$ to be the zero matrix. 
We call $A_i$ the $i$th {\it distance matrix} of $\G$. We abbreviate $A:=A_1$ and call this the 
{\it adjacency matrix} of $\G.$  
We observe
(ai)   $A_0 = I$;
(aii)  $\sum_{i=0}^D A_i = J$;
(aiii) $\overline{A_i} = A_i \;(0 \le i \le D)$;
(aiv)  $A_i^t = A_i  \;(0 \le i \le D)$;
(av)   $A_iA_j = \sum_{h=0}^D p_{ij}^h A_h \;(0 \le i,j \le D)$,
where $I$ (resp. $J$) denotes the identity matrix (resp. all 1's matrix) in  $\MX$. Using these facts
we find $A_0,A_1,\ldots,A_D$ is a basis for a commutative subalgebra $M$ of $\MX$. We call $M$ the 
{\it Bose-Mesner algebra} of $\G$. It turns out that $A$ generates $M$ \cite[p.~190]{BI}. 

\medskip \noindent 
We now recall the dual idempotents of $\G$. To do this
fix a vertex $x \in X.$ We view $x$ as a ``base vertex." For $ 0 \le i \le D$ let 
$E_i^*=E_i^*(x)$ denote the diagonal matrix in $\MX$ with $(y,y)$-entry 
\begin{eqnarray*}
\label{den0}
(\Es_i)_{y y} = \left\{ \begin{array}{lll}
                 1 & \hbox{if } \; \partial(x,y)=i, \\
                 0 & \hbox{if } \; \partial(x,y) \ne i \end{array} \right. 
                 \qquad (y \in X).
\end{eqnarray*}
We call $E_i^*$ the  $i$th 
{\it dual idempotent} of $\G$ with respect to $x$ \cite[p.~378]{ter1}. We observe
(ei)   $\sum_{i=0}^D E_i^*=I$;
(eii)  $\overline{E_i^*} = E_i^*$ $(0 \le i \le D)$;
(eiii) $E_i^{*t} = E_i^*$ $(0 \le i \le D)$;
(eiv)  $E_i^*E_j^* = \delta_{ij}E_i^* $ $(0 \le i,j \le D)$.
By these facts $E_0^*,E_1^*, \ldots, E_D^*$ form a basis for a commutative subalgebra $M^*=M^*(x)$ 
of $\MX$. We call $M^*$ the {\it dual Bose-Mesner algebra} of $\G$ with respect to 
$x$ \cite[p.~378]{ter1}. For $0 \le i \le D$ we have
$$
\Es_i V = {\rm Span} \{ \h{y} \mid y \in X, \: \partial(x,y)=i\},
$$
so ${\rm dim} \Es_i V = k_i$. We call $\Es_i V$ the $i$th {\em subconstituent of} $\G$ 
with respect to $x$. Note that
\begin{equation}
\label{vsub}
V = E_0^*V + E_1^*V + \cdots + E_D^*V \qquad \qquad {\rm (orthogonal\ direct\ sum}).
\end{equation}
Moreover $\Es_i$ is the projection from $V$ onto $\Es_i V$ for $0 \le i \le D$.

\medskip \noindent
We recall the Terwilliger algebra of $\G$. Let $T=T(x)$ denote the subalgebra of $\MX$ generated by 
$M$, $M^*$. We call $T$ the {\it Terwilliger algebra of} $\G$ with respect to $x$ 
\cite[Definition 3.3]{ter1}. Recall $M$ is generated by $A$ so $T$ is 
generated by $A$ and the dual idempotents. We observe $T$ has finite dimension. By construction $T$ is closed under the 
conjugate-transpose map so $T$ is semisimple \cite[Lemma 3.4(i)]{ter1}. 

\medskip \noindent
By a $T$-{\em module} we mean a subspace $W$ of $V$ such that $BW \subseteq W$ 
for all $B \in T$. Let $W$ denote a $T$-module. Then $W$ is said to be {\em irreducible} whenever
$W$ is nonzero and $W$ contains no $T$-modules other than $0$ and $W$.

By \cite[Corollary 6.2]{Go1} any $T$-module is an orthogonal direct sum of irreducible $T$-modules. 
In particular the standard module $V$ is an orthogonal direct sum of irreducible $T$-modules. Let 
$W, W'$ denote $T$-modules. By an {\em isomorphism of} $T$-{\em modules} from $W$ to $W'$ we mean an 
isomorphism of vector spaces $\sigma: W \to W'$ such that $(\sigma B - B \sigma) W = 0$ for all 
$B \in T$. The $T$-modules $W$, $W'$ are said to be {\em isomorphic} whenever there exists an 
isomorphism of $T$-modules from $W$ to $W'$. By \cite[Lemma 3.3]{Cu1} any two nonisomorphic 
irreducible $T$-modules are orthogonal. Let $W$ denote an irreducible $T$-module. By 
\cite[Lemma 3.4(iii)]{ter1} $W$ is an orthogonal direct sum of the nonvanishing spaces among 
$\Es_0 W, \Es_1 W, \ldots , \Es_D W$. By the {\em endpoint} of $W$ we mean 
$\min \{i \mid 0 \le i\le D, \; \Es_i W \ne 0 \}$. By the {\em diameter} of $W$ we mean 
$|\{i \mid 0 \le i \le D, \; \Es_i W \ne 0\}| -1$.   We say $W$ 
is {\it thin} whenever the dimension of $E^*_iW$ is at most 1 for $0 \leq i \leq D$.

Fix a decomposition of the standard module $V$ into an orthogonal direct sum of irreducible $T$-modules.  For any irreducible $T$-module $W$, the {\em multiplicity} of $W$ is the number of irreducible modules in this decomposition which are isomorphic to $W$.  It is well-known that the multiplicity is independent of the decomposition of $V$.  

By \cite[Proposition 8.4]{egge1} $\G$ has a unique irreducible $T$-module with 
endpoint 0. We denote this $T$-module by $V_0$. We call $V_0$ the 
{\em primary module}. It appears in $V$ with multiplicity $1$ and it has basis 
$\{s_i \mid 0 \le i \le D\}$, where 
\begin{equation}
\label{basis}
  s_i = \sum_{y \in \G_i(x)}\hat{y}.
\end{equation}

Recall $\G$ is {\em bipartite} whenever
$a_i=0$ for $0 \le i \le D$. For the rest of this paper we assume $\G$ is bipartite.
In order to avoid trivialities we assume the valency $k \ge 3$.
We now recall a basic formula. Setting $a_i=0$ in (\ref{kk}) we find
\begin{equation}
\label{k_bip}
b_i + c_i = k \quad (0 \le i \le D).
\end{equation}

In the rest of the paper we will consider the following situation.
\begin{notation}
\label{blank} 
Let $\G=(X,R)$ denote a bipartite distance-regular graph with diameter $D \ge 4$, 
valency $k \ge 3$, intersection numbers $b_i, c_i$, and distance matrices $A_i \; (0 \le i \le D)$. 
Let $V$ denote the standard module for $\G$.
We fix $x \in X$ and let $\Es_i = \Es_i(x) \; (0 \le i \le D)$ and $T=T(x)$
denote the dual idempotents and the Terwilliger algebra of $\G$ with respect to $x$, respectively.
Let $V_0$ denote the irreducible $T$-module with endpoint $0$, and let $V_1$ denote the subspace of $V$ spanned by the irreducible 
$T$-modules with endpoint $1$. We define the set $U$ to be the orthogonal complement of 
$\Es_2 V_0 + \Es_2 V_1$ in $\Es_2 V$. 
\end{notation}

With reference to Notation \ref{blank} we now define a certain
partition of $X$ that we will find useful.

\begin{definition}
\label{def:D}
With reference to Notation \ref{blank} fix vertex $y \in X$ such that $\partial(x,y) = 2$.
For all integers $i,j$ define $D_j^i=D_j^i(x,y)$ by
$$
  D_j^i = \{z \in X \mid \partial(x,z)=i \; \hbox{ and } \; \partial(y,z)=j\}.
$$
\end{definition}
We observe $D_j^i = \emptyset$ unless $0 \le i,j \le D$ and $i+j$ is even.


\section{Local eigenvalues}
\label{loceig}

Later in the paper we will consider the irreducible $T$-modules with endpoint $2$.
In order to discuss these we introduce some parameters we call local eigenvalues.

\begin{definition}
\label{def:G22}
With reference to Notation \ref{blank} let $\G_2^2 = \G_2^2(x)$ denote the graph with vertex set 
$\tilde{X}=\G_2(x)$ and edge set $\tilde{R} = \{yz \mid y, z \in \tilde{X}, \; \partial(y,z)=2\}$.
The graph $\G_2^2$ has exactly $k_2$ vertices and it is regular with valency $p_{22}^2$. 
Let $\tilde{A}$ denote the adjacency matrix of $\G_2^2$. The matrix $\tilde{A}$ is symmetric with
real entries.  Therefore $\tilde{A}$ is diagonalizable with all eigenvalues real. Let 
$\eta_1, \eta_2, \ldots, \eta_{k_2}$ denote the eigenvalues of $\tilde{A}$. We call 
$\eta_1, \eta_2, \ldots, \eta_{k_2}$ the {\em local eigenvalues of} $\G$ {\em with respect to} $x$. 
\end{definition}

With reference to Notation \ref{blank} consider the second subconstituent $\Es_2 V$. We recall the 
dimension of $\Es_2 V$ is $k_2$ and observe $\Es_2 V$ is invariant under multiplication by $\Es_2 A_2 \Es_2$.
Note that for an appropriate ordering of the vertices of $\G$ we have 
$$
  \Es_2 A_2 \Es_2 = \begin{pmatrix} \tilde{A} & 0 \\0 & 0 \end{pmatrix},
$$
where $\tilde{A}$ is as in Definition \ref{def:G22}. Apparently the action of $\Es_2 A_2 \Es_2$
on $\Es_2 V$ is essentially the adjacency map for $\G_2^2$. In particular, the action of 
$\Es_2 A_2 \Es_2$ on $\Es_2 V$ is diagonalizable with eigenvalues $\eta_1, \eta_2, \ldots, \eta_{k_2}$.
Note that the vector $s_2$ from \eqref{basis}
is in $\Es_2 V$ and it is an eigenvector for $\Es_2 A_2 \Es_2$ with eigenvalue $p_{22}^2$. 
By \cite[Corollary 7.7]{Cu1} there exists, up to isomorphism, a unique irreducible $T$-module
with endpoint $1$. It has diameter $D-2$ and it appears in $V$ with multiplicity $k-1$. Let $W_1$ denote an irreducible 
$T$-module with endpoint $1$ and choose nonzero $v \in \Es_2 W_1$. By \cite[Corollary 9.1]{MT}, 
$v$ is an eigenvector for $\Es_2 A_2 \Es_2$ with eigenvalue $b_3-1$. Reordering the local eigenvalues 
if necessary, we have $\eta_1 = p_{22}^2$ and $\eta_i = b_3 -1 \; (2 \le i \le k)$.  We define the set $\Phi_2 := \{ \eta_i \mid k+1 \le i \le k_2\}.$
Note that the set $\Phi_2$ is the same as the set $\Phi_2$ from \cite[Definition 4.8]{Cu3}.

\begin{lemma}
\label{lem:loc1}
With reference to Notation \ref{blank}, the local eigenvalues $\eta \in \Phi_2$
are exactly the eigenvalues of $\Es_2 A_2 \Es_2$ on $U$.
\end{lemma}
\proof
By the definition of $U$, we find that $U$ is invariant under the action of $\Es_2 A_2 \Es_2$.
Apparently the restriction of $\Es_2 A_2 \Es_2$ to $U$ is diagonalizable with eigenvalues 
$\eta_i \; (k+1 \le i \le k_2)$. The result follows. \qed


%

\bigskip \noindent
We now define certain scalars $\Delta_i \; (2 \le i \le D-1)$.
To do this recall that, by \cite[Lemma 4.1.7]{BCN}, we have 
$$
  p^i_{2i} = \frac{c_i(b_{i-1}-1) + b_i(c_{i+1}-1) }{c_2}
$$
for $1 \le i \le D-1$.

\begin{definition}
\label{def:delta}
With reference to Notation \ref{blank}, for $2 \le i \le D-1$ we define
$$
  \Delta_i = (b_{i-1} - 1)(c_{i+1} - 1) - (c_2 - 1)p^i_{2i}.
$$
\end{definition}

\begin{lemma}  {\rm (\cite[Theorem 12]{Cu0}, \cite[Corollary 4.13]{Cu3})}
\label{Del2}
With reference to Notation \ref{blank} and Definition \ref{def:delta}, 
we have $\Delta_i \ge 0$ for $2 \le i \le D-1$.
Moreover $\Delta_2 =0$ if and only if $|\Phi_2| \le 1$.
\end{lemma}

\medskip \noindent

Let $W$ denote a thin irreducible $T$-module with endpoint $2$. Then $\Es_2 W$ is a one-dimensional 
eigenspace for $\Es_2 A_2 \Es_2$; we call the corresponding eigenvalue the {\em local eigenvalue}
of $W$. Note this local eigenvalue of $W$ is contained in the set 
$\{\eta_{k+1}, \eta_{k+2}, \ldots, \eta_{k_2}\}$. We will need the following result.

\begin{lemma}  {\rm (\cite[Lemma 10.10 and Theorem 11.7]{Cu2})}
\label{isomequiv}  With reference to Notation \ref{blank},
let $W, W'$ denote thin irreducible $T$-modules with endpoint $2$. 
Then $W$ and $W'$ are isomorphic as $T$-modules if and only if 
they have the same local eigenvalue.
\end{lemma}


\section{Lowering and raising matrices}

With reference to Notation \ref{blank}, in this section we recall the lowering matrix and the 
raising matrix of the algebra $T$.
 
\begin{definition}
\label{def:LR}
{\rm With reference to Notation \ref{blank} we define matrices $L=L(x)$, $R=R(x)$ by
$$
  L = \sum_{h=1}^D \Es_{h-1} A \Es_h, \qquad \qquad
  R = \sum_{h=0}^{D-1} \Es_{h+1} A \Es_h.
$$
Note that $A=L+R$ \cite[Lemma 4.4]{Cu1}. We call $L$ and $R$ the {\em lowering matrix}
and the {\em raising matrix of} $\G$ with respect to $x$, respectively.}
\end{definition}

\begin{lemma}
\label{LR:lem1}
With reference to Notation \ref{blank} let $y,z \in X$. 
Then the following {\rm (i), (ii)} hold.
\begin{itemize}
\item[{\rm (i)}]
$L_{zy} = 1$ if $\partial(z,y)=1$ and $\partial(x,z)=\partial(x,y)-1$, and $0$ otherwise.
\item[{\rm (ii)}]
$R_{zy} = 1$ if $\partial(z,y)=1$ and $\partial(x,z)=\partial(x,y)+1$, and $0$ otherwise.
\end{itemize}
\end{lemma}
\proof Immediate from Definition \ref{def:LR} and elementary matrix multiplication. \qed

\medskip \noindent
With reference to Notation \ref{blank} let $y,z \in X$. We display the $(z,y)$-entry of 
certain products of the matrices $L$ and $R$. To do this we need another definition.

\medskip \noindent
A sequence of vertices $[y_0, y_1, \ldots, y_t]$ of $\Gamma$ is a {\em walk} in $\G$ if $y_{i-1} y_i$ is an edge
for $1 \le i \le t$. 

\begin{lemma}
\label{LR:lem2}
With reference to Notation \ref{blank} choose $y,z \in X$ and let $m$ denote a positive integer.
Assume that $y \in \G_i(x)$. Then the following {\rm (i)--(iii)} hold.
\begin{itemize}
\item[{\rm (i)}]
The $(z,y)$-entry of $R^m$ is equal to the number of walks $[y=y_0, y_1, \ldots, y_m = z]$, 
such that $y_j \in \G_{i+j}(x)$ for $0 \le j \le m$.
\item[{\rm (ii)}]
The $(z,y)$-entry of $R^mL$ is equal to the number of walks $[y=y_0, y_1, \ldots, y_{m+1} = z]$, 
such that $y_1 \in \G_{i-1}(x)$ and $y_j \in \G_{i-2+j}(x)$ for $2 \le j \le m+1$.
\item[{\rm (iii)}]
The $(z,y)$-entry of $LR^m$ is equal to the number of walks $[y=y_0, y_1, \ldots, y_{m+1} = z]$, 
such that $y_j \in \G_{i+j}(x)$ for $0 \le j \le m$ and $y_{m+1} \in \G_{i+m-1}(x)$.
\end{itemize}
\end{lemma}
\proof
Immediate from Lemma \ref{LR:lem1} and elementary matrix multiplication. \qed 

\begin{lemma} \label{LR:zy-entry}
With reference to Notation \ref{blank}, choose an integer $i$ $(2 \le i \le D-2)$. Let $y, z \in X$ such that $y \in \G_2(x), z \in \G_i(x)$.  Then
the following {\rm (i)-(iii)} hold:
\begin{itemize}
\item[{\rm (i)}]
$\displaystyle{ (LR^{i-1})_{z y} = \left\{ \begin{array}{ll}
                 b_i c_{i-1}c_{i-2} \cdots c_1 & \hbox{ if } \; \partial(z,y)=i-2, \\
                 \left(c_i - |\G_{i-1}(x) \cap \G_{i-1}(y) \cap \G(z)|\right)c_{i-1}c_{i-2}\cdots c_1 & \hbox{ if } \; \partial(z,y)=i, \\
                 0 & \hbox{ otherwise. } \end{array} \right.}$
 \item[{\rm (ii)}]
$\displaystyle{ (R^{i-1}L)_{z y} = \left\{ \begin{array}{ll}
                 c_2 c_{i-1}c_{i-2} \cdots c_1 & \hbox{ if } \; \partial(z,y)=i-2, \\
                 |\G(x) \cap \G(y) \cap \G_{i-1}(z)| c_{i-1}c_{i-2} \cdots c_1 & \hbox{ if } \; \partial(z,y)=i, \\
                 0 & \hbox{ otherwise. } \end{array} \right.}$
 \item[{\rm (iii)}]
$\displaystyle{ (R^{i-2})_{z y} = \left\{ \begin{array}{ll}
                 c_{i-2}c_{i-3} \cdots c_1 & \hbox{ if } \; \partial(z,y)=i-2, \\
                 0 & \hbox{ otherwise. } \end{array} \right.}$
 \end{itemize}
 \end{lemma}
 \proof  Since $\partial(x,y)=2, \partial(x,z)=i$, and $\G$ is bipartite, we find $\partial(z,y) \in \{i-2, i, i+2\}$.  
 Now by Lemma \ref{LR:lem2}, the $(z,y)$-entry of each of $LR^{i-1}, R^{i-1}L,$ and $R^{i-2}$ is $0$ whenever $\partial(z,y) =i+2$.  Let $D^\ell_j = D^\ell_j(x,y) \; (0 \le \ell, j \le D)$ be as in Definition \ref{def:D}.

Assume $\partial(z,y)=i$.
By Lemma \ref{LR:lem2}, $(L R^{i-1})_{zy}$ equals the number of walks $[y=y_0, y_1, \ldots, y_i=z]$,
such that $y_j \in \G_{2+j}(x)$ for $0 \le j \le i-1$. Since $\partial(z,y)=i$ each 
such walk is actually a path passing through $D^{i+1}_{i-1}$. Note that $z$ has
$c_i - |\G_{i-1}(x) \cap \G_{i-1}(y) \cap \G(z)|$ neighbours in $D^{i+1}_{i-1}$, and that there are
precisely $c_{i-1} c_{i-2} \cdots c_1$ paths of length $i-1$ from each vertex in $D^{i+1}_{i-1}$
to $y$. It follows from the above comments that
$$
  ( L R^{i-1} )_{zy} = 
  (c_i - |\G_{i-1}(x) \cap \G_{i-1}(y) \cap \G(z)|) c_{i-1} c_{i-2} \cdots c_1.
$$

\noindent
By Lemma \ref{LR:lem2} $(R^{i-1} L)_{zy}$ equals the number of walks $[y=y_0, y_1, \ldots, y_i=z]$,
such that $y_j \in \G_{j}(x)$ for $1 \le j \le i$. Since $\partial(z,y)=i$ each 
such walk is actually a path passing through $D^1_1$. Note that there are
precisely $c_{i-1} c_{i-2} \cdots c_1$ paths of length $i-1$ from $z$ to each vertex in $D^1_1$
which is at distance $i-1$ from $z$. It follows from the above comments that
$$
  ( R^{i-1} L )_{zy} = 
  |\G(x) \cap \G(y) \cap \G_{i-1}(z)| c_{i-1} c_{i-2} \cdots c_1.
$$

\noindent
By Lemma \ref{LR:lem2} $(R^{i-2})_{zy}$ equals the number of walks $[y=y_0, y_1, \ldots, y_{i-2}=z]$,
such that $y_j \in \G_{2+j}(x)$ for $0 \le j \le i-2$. Since $\partial(z,y)=i$, this number equals 0.

Assume finally that $\partial(z,y)=i-2$. Using 
Lemma \ref{LR:lem2} and similar reasoning as above, we find that
$(L R^{i-1})_{zy} = b_i c_{i-1} c_{i-2} \cdots c_1$,
$(R^{i-1} L)_{zy} = c_2 c_{i-1}c_{i-2} \cdots c_1$, and 
$(R^{i-2})_{zy} = c_{i-2} c_{i-3} \cdots c_1$. This completes the proof. \qed

\begin{lemma}
\label{LR:lem4}
With reference to Notation \ref{blank}, the following {\rm (i), (ii)} hold for $2 \le i \le D$:
\begin{itemize}
\item[{\rm (i)}]
$\displaystyle{
  \Es_i A_{i-2}\Es_2 = \frac{1 }{c_1 c_2 \cdots c_{i-2}} \Es_i R^{i-2}\Es_2.}$
\item[{\rm (ii)}]
$\displaystyle{
  \Es_i A_i \Es_2 = 
  \frac{1 }{ c_1 c_2 \cdots c_i} \Es_i \Big(
  R^{i-1} L  + R^{i-2} L R  + \cdots + L R^{i-1}   -
  \sum_{j=1}^{i-1} b_{j-1} c_j  R^{i-2}  \Big)\Es_2.  }$
\end{itemize}
\end{lemma}
\proof
It is well-known that $A_i = v_i(A)$, where $v_j \;(0  \le j \le D)$ are polynomials of degree $j$
defined recursively by
\begin{equation}
\label{poli}
  v_0(\lambda) = 1, \qquad v_1(\lambda) = \lambda, \qquad
  c_{j+1} v_{j+1}(\lambda) = \lambda v_j(\lambda) - b_{j-1} v_{j-1}(\lambda).
\end{equation}
From \eqref{poli} it is easy to see that for $0 \le j \le D$ we have
$$
  v_j(\lambda) = \frac{1 }{c_1 c_2 \cdots c_j} \lambda^j - 
  \sum_{\ell=1}^{j-1} \frac{b_{\ell-1} c_\ell }{ c_1 c_2 \cdots c_j} \lambda^{j-2} + (\mbox{lower degree terms}).
$$
Since $A=R+L$ it follows from the above comments that
$$
A_{i-2} = v_{i-2}(R+L) = \frac{1 }{ c_1 c_2 \cdots c_{i-2}} (R+L)^{i-2} - 
    \sum_{\ell=1}^{i-3} \frac{b_{\ell-1} c_\ell }{ c_1 c_2 \cdots c_{i-2}} (R+L)^{i-4} + \cdots
$$
and
$$
  A_i = v_i(R+L) = \frac{1 }{ c_1 c_2 \cdots c_i} (R+L)^i - 
    \sum_{\ell=1}^{i-1} \frac{b_{\ell-1} c_\ell }{ c_1 c_2 \cdots c_i} (R+L)^{i-2} + \cdots
$$
The result now follows from the facts that $R \Es_j V \subseteq \Es_{j+1} V$ and 
$L \Es_j V \subseteq \Es_{j-1} V$. \qed

\begin{lemma}
\label{lem:loc2}
With reference to Notation \ref{blank}, pick $v \in U$.
Then the following {\rm (i), (ii)} hold.
\begin{itemize}
\item[{\rm (i)}]
The vector $v$ is an eigenvector of $\Es_2 A_2 \Es_2$ if and only if $v$ is an eigenvector of $\Es_2 L R \Es_2$.
\item[{\rm (ii)}]
If $v$ is an eigenvector of $\Es_2 A_2 \Es_2$ with eigenvalue $\eta \in \Phi_2$,
then the corresponding eigenvalue for $\Es_2 LR \Es_2$ is $c_2 \eta + k$.
\end{itemize}
\end{lemma}
\proof
(i) By Lemma \ref{LR:lem4}(ii) we have
$$
  \Es_2 A_2 \Es_2 = \frac{1 }{c_2} \Es_2 (RL + LR - kI) \Es_2.
$$
Moreover, since $v \in U$, we have $\Es_2 RL \Es_2 v = RL v = 0$. The result follows. 

\smallskip \noindent
(ii) Immediate from comments in the proof of (i) above. \qed


\section{Algebraic condition implies combinatorial property}

With reference to Notation \ref{blank} assume that up to isomorphism $\G$ has exactly
two irreducible $T$-modules with endpoint $2$, and they are both thin. In this section
we prove that for $2 \le i \le D - 2$ there exist complex scalars $\alpha_i, \beta_i,\gamma_i$, 
not all zero, with the following property: for all $y, z \in X$ such that 
$\partial(x, y) = 2, \: \partial(x, z) = i, \: \partial(y, z) = i$ we have
$$
  \alpha_i + \beta_i  |\G(x) \cap \G(y) \cap \G_{i-1}(z)| + 
  \gamma_i |\G_{i-1}(x) \cap \G_{i-1}(y) \cap \G(z)| = 0.
$$

Since up to isomorphism $\G$ has exactly two irreducible $T$-modules with endpoint $2$,
we can split the irreducible $T$-modules with endpoint $2$ into two isomorphism classes.
Let $V_2^1$ denote the sum of $T$-modules from the first of these isomorphism classes, 
and let $V_2^2$ denote the sum of $T$-modules from the second of these isomorphism classes.
Let $d_1$ denote the diameter of the irreducible $T$-modules from the first isomorphism class,
and let $d_2$ denote the diameter of the irreducible $T$-modules from the second isomorphism class.
Note that $d_1, d_2 \in \{D-4,D-3,D-2\}$ by \cite[Theorem 10.1]{Cu2}.

\begin{lemma}
\label{partI:lem0}
With reference to Notation \ref{blank} assume that up to isomorphism $\G$ has exactly 
two irreducible $T$-modules with endpoint $2$, and they are both thin. Let 
$F_1$, $F_2$, $F_3$, $F_4$, $F_5 \in T$ and pick an integer $i, \; 2 \le i \le D-2$. 
Assume that $\Es_i F_5 \Es_2 v$ is nonzero for every nonzero $v \in \Es_2 V$.
Then the matrices 
$$
  \Es_i F_1 \Es_2, \; \Es_i F_2 \Es_2, \; \Es_i F_3 \Es_2, \; 
  \Es_i F_4 \Es_2, \Es_i F_5 \Es_2
$$
are linearly dependent.
\end{lemma}
\proof
Since $V_0$ is thin, there exist scalars $r_{0,j} \; (1 \le j \le 4)$ such that for 
$v \in \Es_2 V_0$ we have $\Es_i F_j \Es_2 v = r_{0,j} \Es_i F_5 \Es_2 v$. 
Since all irreducible $T$-modules with endpoint $1$ are thin and mutually isomorphic,
there exist scalars $r_{1,j} \; (1 \le j \le 4)$ such that for 
$v \in \Es_2 V_1$ we have $\Es_i F_j \Es_2 v = r_{1,j} \Es_i F_5 \Es_2 v$. 
Since all the summands of $V_2^\ell \; (\ell \in \{1,2\})$ are thin and mutually isomorphic,
there exist scalars $r_{2,j}^\ell \; (1 \le j \le 4)$ such that for 
$v \in \Es_2 V_2^\ell$ we have $\Es_i F_j \Es_2 v = r_{2,j}^\ell \Es_i F_5 \Es_2 v$. 

We will now show that there exist scalars $\lambda_i \; (1 \le i \le 5)$ such that $\lambda_i$ are not
all zero and such that
\begin{equation}
\label{dep}
  \lambda_1 \Es_i F_1 \Es_2 + \lambda_2 \Es_i F_2 \Es_2 + \lambda_3 \Es_i F_3 \Es_2 + 
  \lambda_4 \Es_i F_4 \Es_2 + \lambda_5 \Es_i F_5 \Es_2 = 0.
\end{equation}
Let $S$ denote the matrix 
$$
\left( \begin{array}{ccccc}
	    r_{0,1}   & r_{0,2}   & r_{0,3}   & r_{0,4}   & 1 \cr
            r_{1,1}   & r_{1,2}   & r_{1,3}   & r_{1,4}   & 1 \cr
            r_{2,1}^1 & r_{2,2}^1 & r_{2,3}^1 & r_{2,4}^1 & 1 \cr
            r_{2,1}^2 & r_{2,2}^2 & r_{2,3}^2 & r_{2,4}^2 & 1 \cr
          \end{array}\right).
$$
Since the rank of $S$ is less or equal to $4$, the system of equations 
$S (x_1, x_2, x_3, x_4, x_5)^t = (0,0,0,0,0)^t$ has a nontrivial solution.
Let $(\lambda_1, \lambda_2, \lambda_3, \lambda_4, \lambda_5)^t$ be such a nontrivial solution.
Choose an arbitrary $v \in V$ and note that 
$$
  \Es_2 v = v_0 + v_1 + v_2^1 + v_2^2,
$$
where $v_0 \in V_0$, $v_1 \in V_1$, $v_2^1 \in V_2^1$ and $v_2^2 \in V_2^2$.
Using the above comments we find
$$
  \Es_i F_j \Es_2 v = r_{0,j} \Es_i F_5 \Es_2 v_0 + r_{1,j} \Es_i F_5 \Es_2 v_1 +
                      r_{2,j}^1 \Es_i F_5 \Es_2 v_2^1 + r_{2,j}^2 \Es_i F_5 \Es_2 v_2^2
                      \qquad (1 \le j \le 4).
$$
It follows that
$$
  \lambda_1 \Es_i F_1 \Es_2 v + \lambda_2 \Es_i F_2 \Es_2 v + \lambda_3 \Es_i F_3 \Es_2 v + 
  \lambda_4 \Es_i F_4 \Es_2 v + \lambda_5 \Es_i F_5 \Es_2 v = 0
$$
and hence \eqref{dep} holds. \qed

\begin{lemma}
\label{partI:lem2}
With reference to Notation \ref{blank} assume that up to isomorphism $\G$ has exactly 
two irreducible $T$-modules with endpoint $2$, and they are both thin. For $2 \le i \le D-2$,
the matrices 
$$
  \Es_i L R^{i-1} \Es_2, \; \Es_i R^{i-1} L \Es_2, \;  \Es_i A_i \Es_2, \; \Es_i R^{i-2} \Es_2
$$
are linearly dependent.
\end{lemma}
\proof
By \cite[Lemma 5.3(iii), Lemma 7.3(ii)]{Cu1} and \cite[Theorem 9.4(i), Lemma 10.2(i)]{Cu2}, 
$R^{i-2}v \ne 0$ for each nonzero $v \in \Es_2 V$. Therefore, by Lemma \ref{partI:lem0},
there exist scalars $\lambda_j \; (1 \le j \le 5)$, not all zero,
such that 
$$
  \lambda_1 \Es_i R^i L^2 \Es_2  + \lambda_2 \Es_i L R^{i-1} \Es_2  + \lambda_3  \Es_i R^{i-1} L \Es_2+
  \lambda_4  \Es_i A_i \Es_2 + \lambda_5  \Es_i R^{i-2} \Es_2= 0.
$$
We will show that $\lambda_1=0$.
Choose $y,z \in X$ such that $\partial(x,y)=2$, $\partial(x,z)=i$ and $\partial(y,z)=i+2$.
Then by Lemma \ref{LR:lem2} and elementary matrix multiplication, the $(z,y)$-entry of $\Es_i R^i L^2 \Es_2$ is nonzero, while by Lemma \ref{LR:zy-entry} the 
$(z,y)$-entries of 
$\Es_i L R^{i-1} \Es_2$, $\Es_i R^{i-1} L \Es_2$, and $\Es_i R^{i-2} \Es_2$ 
are all 0. Similarly, the  $(z,y)$-entry of 
$\Es_i A_i \Es_2$ is 0.  Hence $\lambda_1=0$ and the result follows. \qed

\begin{theorem}
\label{partI:thm_mainI}
With reference to Notation \ref{blank} assume that up to isomorphism $\G$ has exactly 
two irreducible $T$-modules with endpoint $2$, and they are both thin. Then
for $2 \le i \le D - 2$ there exist complex scalars $\alpha_i, \beta_i, \gamma_i$ which are not all zero 
and with the following property: for all $y, z \in X$ such that 
$\partial(x, y) = 2, \: \partial(x, z) = i, \: \partial(y, z) = i$ we have
\begin{equation} \label{eq:recurrence}
  \alpha_i + \beta_i  |\G(x) \cap \G(y) \cap \G_{i-1}(z)| + 
  \gamma_i |\G_{i-1}(x) \cap \G_{i-1}(y) \cap \G(z)| = 0.
\end{equation}
\end{theorem}
\proof
Note that the theorem holds for $i=2$ with $\alpha_2=0$, $\beta_2=1$ and $\gamma_2=-1$. Assume $i \ge 3$.
By Lemma \ref{partI:lem2} there exist scalars $\lambda_{1,i}, \lambda_{2,i}, \lambda_{3,i}$ and $\lambda_{4,i}$ 
which are not all zero and such that
\begin{equation}
\label{eq:dep}
  \lambda_{1,i} \Es_i L R^{i-1} \Es_2 + \lambda_{2,i} \Es_i R^{i-1} L \Es_2 + 
  \lambda_{3,i} \Es_i A_i \Es_2 + \lambda_{4,i}  \Es_i R^{i-2} \Es_2= 0.
\end{equation}
Let $y,z \in X$ such that $\partial(x,y)=2$, $\partial(x,z)=i$ and $\partial(y,z)=i$.
Observe that since $\partial(x,y)=2$ and $\partial(x,z)=i$ we have 
$(\Es_i L R^{i-1} \Es_2)_{zy} = (L R^{i-1})_{zy}$, $(\Es_i R^{i-1} L \Es_2)_{zy} = (R^{i-1} L)_{zy}$,
$(\Es_i A_i \Es_2)_{zy} = (A_i)_{zy}$, and $(\Es_i R^{i-2} \Es_2)_{zy} = (R^{i-2})_{zy}$.
Now observe that $(\Es_i A_i \Es_2)_{zy}=1$ by (\ref{dm1}).
By Lemma \ref{LR:zy-entry}, we find
\begin{eqnarray*}
  (\Es_i L R^{i-1} \Es_2)_{zy} &=&
  (c_i - |\G_{i-1}(x) \cap \G_{i-1}(y) \cap \G(z)|) c_{i-1} c_{i-2} \cdots c_1, \\
  (\Es_i R^{i-1} L \Es_2)_{zy} &=&
  |\G(x) \cap \G(y) \cap \G_{i-1}(z)| c_{i-1} c_{i-2} \cdots c_1, \\
  (\Es_i R^{i-2} \Es_2)_{zy} &=& 0.
\end{eqnarray*}

Set $\alpha_i =  \lambda_{1,i} c_i \cdots c_1 + \lambda_{3,i}$, 
$\beta_i = \lambda_{2,i} c_{i-1} \cdots c_1$ and
$\gamma_i = -\lambda_{1,i} c_{i-1} \cdots c_1$. 
If $\alpha_i = \beta_i = \gamma_i = 0$, then $\lambda_{1,i} = \lambda_{2,i} = \lambda_{3,i} = 0$.
It follows from \eqref{eq:dep} that $\lambda_{4,i} \Es_i R^{i-2} \Es_2 = 0$. Since 
$\lambda_{1,i}, \lambda_{2,i}, \lambda_{3,i}, \lambda_{4,i}$ are not all zero, 
we have $\Es_i R^{i-2} \Es_2=0$. But this is 
a contradiction, since by Lemma \ref{LR:lem2}(i), $(\Es_i R^{i-2} \Es_2)_{wy}\ne 0$ for $w \in \G_i(x) \cap \G_{i-2}(y)$. 
Therefore $\alpha_i, \beta_i$ and $\gamma_i$ are not all zero.  Taking the $(z,y)$-entry of both sides of (\ref{eq:dep}) and using the above information, we obtain the desired result.  
\qed

\section{The $\gamma_i=0$ case}

With reference to Notation \ref{blank}, assume that up to isomorphism $\G$ has exactly 
two irreducible $T$-modules with endpoint $2$, and they are both thin. In the last section we showed that 
for $2 \le i \le D - 2$ there exist complex scalars $\alpha_i, \beta_i, \gamma_i$ which are not all zero 
and with the following property: for all $y, z \in X$ such that 
$\partial(x, y) = 2, \: \partial(x, z) = i, \: \partial(y, z) = i$ we have
$$
  \alpha_i + \beta_i  |\G(x) \cap \G(y) \cap \G_{i-1}(z)| + 
  \gamma_i |\G_{i-1}(x) \cap \G_{i-1}(y) \cap \G(z)| = 0.
$$

In this section, we will show that $\gamma_i \ne 0$ for $2 \le i \le D-2$.  We will need the following lemma.

\begin{lemma} \label{curtinrecurrence} {\rm (\cite[Section 10]{Cu2})}
With reference to Notation \ref{blank}, let $W$ denote a thin irreducible $T$-module of endpoint $2$.  Let $v \in \Es_2 W$ be an eigenvector for $\Es_2 A_2 \Es_2$ with eigenvalue $\eta$.  For all integers $2 \le  i \le D-1$, define  vectors $v_i^+ = \Es_i A_{i-2}v, \; v_i^- = \Es_i A_{i+2}v$.   There exist unique real scalars $\varphi_i = \varphi_i(W),$ $\omega_i =\omega_i(W)$ $(2 \le i \le D-2)$ such that 
$$v_i^- = \varphi_i v_i^+, \qquad Lv_{i+1}^+ = \omega_i v_i^+.$$
The scalars $\varphi_i$ and $\omega_i$ satisfy the following recurrence:
\begin{eqnarray} \label{rec:beta2}
\varphi_2 = -(\eta+1), \qquad & & \omega_2 = b_2 -c_2 \varphi_2,
 \\
 \label{rec:gammai}
 \varphi_i = \frac{b_{i+1}}{\omega_{i-1}} \varphi_{i-1}, \qquad & &
  \omega_i = c_{i-2} \varphi_{i-1} + b_i - c_i \varphi_i \qquad (3 \le i \le D-2).
 \end{eqnarray}
 Moreover, the scalars $\omega_2, \omega_3, \ldots, \omega_{D-3}$ are all positive.
 \end{lemma}
 
 For the rest of this section, we will use the following notation.
 
 \begin{notation} \label{blank2} 
 With reference to Notation \ref{blank}, assume that up to isomorphism $\G$ has exactly 
two irreducible $T$-modules with endpoint $2$, and they are both thin.  Let $\alpha_i, \beta_i, \gamma_i \; (2 \le i \le D-2)$ denote the scalars from Theorem \ref{partI:thm_mainI}.  Let $W, W'$ denote non-isomorphic thin irreducible $T$-modules of endpoint $2$ with corresponding local eigenvalues $\eta, \eta'$, respectively.  We define scalars $\psi, \psi'$ by
\begin{equation} \label{psidef}
 \psi = \begin{cases}b_2 - \frac{b_2 b_3}{1+ \eta}, & \mbox{ if }\eta \ne -1 \\
\infty, & \mbox{ if }\eta = -1 \end{cases} \qquad \qquad
\psi' = \begin{cases}b_2 - \frac{b_2 b_3}{1+ \eta'}, & \mbox{ if }\eta' \ne -1 \\
\infty, & \mbox{ if }\eta' = -1. \end{cases}
\end{equation}
Fix nonzero vectors $v \in \Es_2W, v' \in \Es_2W'$.  We note $v, v'$ are eigenvectors for $\Es_2 A_2 \Es_2$ with eigenvalues $\eta, \eta'$, respectively.  We define scalars $\varphi_i(W),  \varphi_i(W'), \omega_i(W), \omega_i(W')$ $(2 \le i \le D-2)$ as in Lemma \ref{curtinrecurrence}.
\end{notation}

With reference to Notation \ref{blank2}, recall that we wish to show $\gamma_i \ne 0$ $(2 \le i \le D-2)$.  We will proceed by contradiction.  As our proof is a bit long, we will need some intermediary results.  In the following two lemmas, we will assume there exists an integer $i$ $(2 \le i \le D-2)$ such that $\gamma_i =0$, and we obtain some preliminary results that ultimately will lead to a contradiction.

\begin{lemma} \label{equalvarphi}
With reference to Notation \ref{blank2}, assume there exists an integer $i$ $(2 \le i \le D-2)$ such that $\gamma_i =0$.   Then 
$\varphi_i(W) = \varphi_i(W')$.  Moreover, these scalars are both nonnegative.  
\end{lemma}
\proof Setting $\gamma_i=0$ in (\ref{eq:recurrence}), we find that $ |\G(x) \cap \G(y) \cap \G_{i-1}(z)|$ is constant for all $y ,z \in X$ such that $\partial(x,y)=2, \partial(x,z)=i, \partial(y,z)=i$.  Let $\tau$ denote this constant.  Now let $z, y$ denote arbitrary vertices in $X$. 
By computing the $(z,y)$-entry of both sides, we shall now show that
\begin{equation} \label{new}
   \Es_i R^{i-1}L\Es_2 - c_{i-1}c_2 \Es_i R^{i-2}\Es_2 -\tau c_{i-1}c_{i-2} \cdots c_1 \Es_i A_i \Es_2 =0.
\end{equation}
Indeed, it is clear that the $(z,y)$-entry of both sides of (\ref{new}) is zero unless 
 $z \in \Gamma_i(x)$ and $y \in \Gamma_2(x)$. So let $z \in \Gamma_i(x)$, $y \in \Gamma_2(x)$.
Observe that 
$(\Es_i R^{i-1} L \Es_2)_{zy} = (R^{i-1} L)_{zy}$, $(\Es_i R^{i-2} \Es_2)_{zy} = (R^{i-2})_{zy}$, and
$(\Es_i A_i \Es_2)_{zy} = (A_i)_{zy}$.
Using Lemma \ref{LR:zy-entry} and (\ref{dm1}), one verifies the $(z,y)$-entries of both sides of (\ref{new}) are equal, and hence (\ref{new}) holds.

Since the $T$-module $W$ has endpoint $2$, we observe $\Es_i R^{i-1}L\Es_2 v=0$ and thus
\begin{equation} \label{eq:varphi1}  
\tau c_{i-1}c_{i-2} \ldots c_1 \Es_i A_i \Es_2 v = - c_{i-1}c_2\Es_i R^{i-2} \Es_2 v.
\end{equation}
By \cite[Corollary 9.3, Theorem 9.4]{Cu2}, we find
\begin{equation} \label{eq:varphi2} 
\Es_i A_i \Es_2 v = -v_i^+-v_i^-, \qquad \Es_i R^{i-2} \Es_2v = c_{i-2}c_{i-3} \ldots c_2 v_i^+,
\end{equation}
where $v_i^+, v_i^-$ are as in Lemma \ref{curtinrecurrence}.
Combining the information in (\ref{eq:varphi1}), (\ref{eq:varphi2}), we find $\tau^{-1}(c_2 - \tau) v_i^+ = v_i^-$.  Thus by Lemma \ref{curtinrecurrence}, $\varphi_i(W) = \tau^{-1}(c_2 - \tau)$.  
We note that $\tau^{-1} (c_2 - \tau) \ge 0$ since $|\Gamma(x) \cap \Gamma(y)| = c_2$.
By a similar argument, $\varphi_i(W') = \tau^{-1}(c_2 - \tau)$, and the result follows.
\qed

\begin{lemma} \label{etanotminus1}
With reference to Notation \ref{blank2}, assume there exists an integer $i$ $(2 \le i \le D-2)$ such that $\gamma_i =0$.   Then $\eta, \eta' \ne -1$.  
\end{lemma}
\proof  Suppose, to the contrary, that $\eta=-1$.  Then by (\ref{rec:beta2}), (\ref{rec:gammai}), we find $\varphi_2(W)=0$, $\varphi_i(W)=0.$  By Lemma \ref{equalvarphi}, we find $\varphi_i(W')=0$, and hence $\varphi_2(W')=0$ by (\ref{rec:gammai}).  So $\eta' = -1$ by (\ref{rec:beta2}).  Thus $W, W'$ have the same local eigenvalue, and hence are isomorphic by Lemma \ref{isomequiv}, a contradiction.  
\qed

\begin{theorem}  
\label{thm:gamma}
With reference to Notation \ref{blank2}, $\gamma_i \ne 0$ $(2 \le i \le D-2)$.
\end{theorem}
\proof  Suppose there exists an integer $i$ such that $\gamma_i=0$.  By Lemma \ref{equalvarphi}, $\varphi_i(W) = \varphi_i(W')$.  By Lemma \ref{etanotminus1}, $\eta, \eta' \ne -1$.  We now define polynomials $p_j, P_j$ $(0 \le j \le D-2)$ in $\RR[\lambda]$ as in \cite[Definition 6.2]{MT2}.  Using 
\cite[Lemma 6.3]{MT2}, (\ref{rec:beta2})--(\ref{psidef}), and treating separately the cases where $j$ is odd and even, one uses induction on $j$ to routinely prove that
\begin{equation} \label{varphi-p_i}
\varphi_j(W) = \frac{b_j b_{j+1}}{c_{j-1}c_j} \frac{P_{j-2}(\psi)}{P_j(\psi)}, \qquad \varphi_j(W') = \frac{b_j b_{j+1}}{c_{j-1}c_j} \frac{P_{j-2}(\psi')}{P_j(\psi')} \qquad (2 \le j \le D-2).
\end{equation}
Now let $\theta_1$ denote the second-largest eigenvalue of $\Gamma$, and let $\theta_d$ denote the smallest nonnegative eigenvalue of $\Gamma$.  In \cite[Theorem 11.4]{MT} it is shown that $\tilde{\theta}_1 \le \eta, \eta' \le \tilde{\theta}_d$, where $\tilde{\theta}_1 = -1-b_2b_3(\theta_1^2-b_2)^{-1}$ and $\tilde{\theta}_d = -1-b_2b_3(\theta_d^2-b_2)^{-1}$.  By \cite[Lemma 3.5]{MT}, $\theta_1^2 > b_2 > \theta_d^2$, so $\tilde{\theta}_1 < -1 < \tilde{\theta}_d$.  We claim that $\tilde{\theta}_1 \le \eta, \eta' < -1.$  To the contrary, without loss of generality, suppose $-1 < \eta \le \tilde{\theta}_d$.  By \cite[Lemmas 5.6(i), 6.6(iii), 15.3(ii)]{MT2}, we find that $P_{i-2}(\psi), P_i(\psi)$ have opposite signs.  Thus in view of (\ref{varphi-p_i}), we find $\varphi_i(W)$ is negative, contradicting Lemma \ref{equalvarphi}.  Thus $\tilde{\theta}_1 \le \eta, \eta' < -1$, and by \cite[Lemmas 5.5(i), 15.3(i)]{MT2}, we find $P_j(\psi), P_j(\psi')$ are positive for $0 \le j \le D-2$.  In view of the fact that $ \eta, \eta' < -1$ and using (\ref{psidef}), we find $\psi, \psi' >0$.

By (\ref{varphi-p_i}) and the fact that $\varphi_i(W) = \varphi_i(W')$, we find
\begin{equation}
 P_{i-2}(\psi) P_i(\psi') -  P_{i-2}(\psi')P_i(\psi)=0.
 \end{equation}
 Now using \cite[Definition 6.2]{MT2}, \cite[Lemma 5.3]{MT}, and treating the cases of $i$ odd and even separately, we find
 \begin{equation}
 c_{i-1}^{-1} c_i^{-1} (\psi' - \psi) \sum_{{h=0}\atop{ i-h {\tiny \mbox{ even}}}}^{i-2} 
P_{h}(\psi) P_{h}(\psi')
\frac{k_{i-2}b_{i-2}b_{i-1}}
{k_{h}b_{h}b_{h+1}}=0.
\end{equation}
Since all terms in the sum above are positive, we find $\psi=\psi'$.  Thus $\eta = \eta'$, so $W, W'$ are isomorphic by Lemma \ref{isomequiv}, a contradiction.
\qed

\begin{corollary}
\label{partI:cor_mainI}
With reference to Notation \ref{blank} assume that up to isomorphism $\G$ has exactly 
two irreducible $T$-modules with endpoint $2$, and they are both thin. Then $\Delta_2 > 0$ and
for $2 \le i \le D - 2$ there exist complex scalars $\alpha_i, \beta_i$ with the following property:  for all $y, z \in X$ such that 
$\partial(x, y) = 2, \: \partial(x, z) = i, \: \partial(y, z) = i$ we have
\begin{equation} \label{eq:recurrence2}
  \alpha_i + \beta_i  |\G(x) \cap \G(y) \cap \G_{i-1}(z)| = 
   |\G_{i-1}(x) \cap \G_{i-1}(y) \cap \G(z)|.
\end{equation}
\end{corollary}

\proof
Since $\G$ has exactly two thin irreducible $T$-modules with endpoint 2, we see $|\Phi_2|=2$ by Lemma \ref{isomequiv}. Thus $\Delta_2 > 0$ by Lemma \ref{Del2}.  Equation (\ref{eq:recurrence2}) follows immediately from Theorem \ref{partI:thm_mainI} and Theorem \ref{thm:gamma}.  
\qed

\section{Some equations involving lowering and raising matrices}

In this section we assume the following notation.

\begin{notation}
\label{blankII}
With reference to Notation \ref{blank} assume that for $2 \le i \le D - 2$ there exist complex scalars 
$\alpha_i, \beta_i$ with the following property: for all $y, z \in X$ such that 
$\partial(x, y) = 2, \: \partial(x, z) = i, \: \partial(y, z) = i$ we have
\begin{equation} \label{main-recur}
  \alpha_i + \beta_i  |\G(x) \cap \G(y) \cap \G_{i-1}(z)| = |\G_{i-1}(x) \cap \G_{i-1}(y) \cap \G(z)|.
\end{equation} 
\end{notation}
Assume the situation from Notation \ref{blankII}.
Our goal is to prove that if $\Delta_2>0$, then up to isomorphism $\G$ has exactly two irreducible $T$-modules 
with endpoint $2$, and they are both thin. To do this, we will need some equations involving raising and lowering matrices.

\begin{lemma}
\label{partII:lem0a}
With reference to Notation \ref{blankII}, for $2 \le i \le D-2$ we have
\begin{equation}
\label{mat_dep}
 L R^{i-1} \Es_2 = c_{i-1} \cdots c_1 (c_i - \alpha_i) \Es_i A_i \Es_2 - 
                        \beta_i  R^{i-1} L \Es_2 +
                        c_{i-1}(b_i + \beta_i c_2)  R^{i-2} \Es_2.
\end{equation}
\end{lemma}
\proof
Let $y, z \in X$. We show that the $(z,y)$-entry of both sides of the above equation agree.
First note that the $(z,y)$-entry of both sides is $0$ if $\partial(x,y) \ne 2$.
It follows from Lemma \ref{LR:lem2} that the $(z,y)$-entry of both sides is $0$ if $\partial(x,z) \ne i$.

Assume now $\partial(x,y) = 2$ and $\partial(x,z) = i$. 
Observe that 
$(LR^{i-1}\Es_2)_{zy} = (LR^{i-1})_{zy}$,  $(\Es_i A_i \Es_2)_{zy} = (A_i)_{zy}$, $(R^{i-1} L \Es_2)_{zy} = (R^{i-1} L)_{zy}$, and $( R^{i-2} \Es_2)_{zy} = (R^{i-2})_{zy}$.
By the triangle inequality and since 
$\G$ is bipartite we have $\partial(y,z) \in \{i-2,i,i+2\}$.
In each of these three cases, we may use Lemma \ref{LR:zy-entry}, (\ref{dm1}), and (\ref{main-recur}) to verify that the $(z,y)$-entries of both sides of \eqref{mat_dep} agree.
This completes the proof. \qed

\begin{corollary}
\label{partII:cor0a}
With reference to Notation \ref{blankII}, for $2 \le i \le D-2$ we have
\begin{equation}
\label{mat_dep1}
\begin{split}
\alpha_i L R^{i-1} \Es_2 & = 
(c_i - \alpha_i) \Big( R^{i-1} L \Es_2 + R^{i-2} L R \Es_2 + \cdots + RLR^{i-2} \Es_2 \Big) + \\
& \quad \;\; -c_i \beta_i R^{i-1} L \Es_2 + 
\Big( c_i c_{i-1} (\beta_i c_2 + b_i) + (\alpha_i - c_i) \sum_{j=1}^{i-1} b_{j-1} c_j \Big) R^{i-2} \Es_2.
\end{split}
\end{equation}
\end{corollary}
\proof
From Lemma \ref{LR:lem4}(ii) we find
\begin{equation}
\label{pomozna1}
\begin{split}
LR^{i-1} \Es_2 & = c_i c_{i-1} \cdots c_1 \Es_i A_i \Es_2 - 
R^{i-1} L \Es_2 - R^{i-2} L R \Es_2 - \cdots - RLR^{i-2} \Es_2 + \\
& \quad \;\; \sum_{j=1}^{i-1} b_{j-1} c_j R^{i-2} \Es_2.
\end{split}
\end{equation}
Multiplying \eqref{pomozna1} with $c_i - \alpha_i$ and \eqref{mat_dep} with $c_i$
and then subtracting the resulting equations, we find \eqref{mat_dep1}. \qed


\section{The $\alpha_i = 0$ case}

With reference to Notation \ref{blankII}, in this section we assume that $\alpha_i=0$
for some $i$ $(3 \le i \le D-2)$. We will show that in this case $\Delta_2 = 0$, where $\Delta_2$ is from Definition \ref{def:delta}.

\begin{definition}
\label{def:blank3}
With reference to Notation \ref{blankII} assume that $\alpha_i=0$ for some $i$ $(3 \le i \le D-2)$, 
and let $\ell \ge 3$ be the minimal integer such that $\alpha_{\ell}=0$. 
For $3 \le i \le D-2$ define $t_i$ as
\begin{equation}
\label{eq:main1}
  t_i = c_i c_{i-1} (\beta_i c_2 + b_i) + (\alpha_i - c_i) \sum_{j=1}^{i-1} b_{j-1} c_j.
\end{equation}
Moreover, let $u \in U$ be an eigenvector for $LR$, and let $\lambda$
be the corresponding eigenvalue. Define scalars $\lambda_i \; (1 \le i \le \ell-2)$ by
\begin{equation}
\label{eq:main2}
  \lambda_1 = \lambda, \qquad \qquad \lambda_i = \frac{c_{i+1} - \alpha_{i+1} }{ \alpha_{i+1}} \Big( \lambda_1 + \lambda_2 + \cdots + \lambda_{i-1} \Big) + 
  \frac{t_{i+1} }{ \alpha_{i+1}}.
\end{equation}
\end{definition}

\medskip

Note that for $1 \le i \le \ell-2$ we have $\lambda_i = \sigma_i \lambda + \rho_i$, where $\sigma_1 = 1, \rho_1 = 0$ and 
\begin{equation}
\label{eq:B}
  \sigma_i = \frac{c_{i+1} - \alpha_{i+1} }{\alpha_{i+1}} \Big( \sigma_1 + \sigma_2 + \cdots + \sigma_{i-1} \Big)
  \qquad \quad (2 \le i \le \ell-2),
\end{equation}
\begin{equation}
\label{eq:A}
  \rho_i = \frac{c_{i+1} - \alpha_{i+1} }{ \alpha_{i+1}} \Big( \rho_1 + \rho_2 + \cdots + \rho_{i-1} \Big) + \frac{t_{i+1} }{\alpha_{i+1}}
  \qquad \quad (2 \le i \le \ell-2).
\end{equation}

We now show that $\sigma_1 + \sigma_2 + \cdots + \sigma_i \ne 0$ for $1 \le i \le \ell-2$.

\begin{lemma}
\label{main:lem6}
With reference to Definition \ref{def:blank3} we have
$\sigma_1 + \sigma_2 + \cdots + \sigma_i \ne 0$ for $1 \le i \le \ell-2$.
\end{lemma}
\proof
When $i = 1$ the result is clear. Assume now that $i \ge 2$ and that on the contrary we have 
$\sigma_1 + \sigma_2 + \cdots + \sigma_i=0$. We claim that in this case also $\sigma_1 + \sigma_2 + \cdots + \sigma_{i-1}=0$.
If $c_{i+1} - \alpha_{i+1} = 0$, then, by \eqref{eq:B}, we have $\sigma_i=0$ and so 
$\sigma_1 + \sigma_2 + \cdots + \sigma_{i-1} = 0$.

If $c_{i+1} - \alpha_{i+1} \ne 0$, then, by \eqref{eq:B}, we have 
$$
  \sigma_1 + \sigma_2 + \cdots + \sigma_{i-1} = \frac{\alpha_{i+1} \sigma_i }{ c_{i+1} - \alpha_{i+1}}.
$$
Therefore
$$
  0 = \sigma_1 + \sigma_2 + \cdots + \sigma_{i-1} + \sigma_i = \frac{\alpha_{i+1} \sigma_i }{c_{i+1} - \alpha_{i+1}} + \sigma_i = 
  \frac{c_{i+1} \sigma_i }{ c_{i+1} - \alpha_{i+1}}.
$$
It follows that $\sigma_i=0$, and so $\sigma_1 + \sigma_2 + \cdots + \sigma_{i-1}=0$.
But now it follows that $\sigma_1=0$, a contradiction.
\qed

\begin{lemma}
\label{main:lem1}
With reference to Definition \ref{def:blank3}, for $1 \le i \le \ell-2$ we have
$$
  LR^i u = \lambda_i R^{i-1} u.
$$
\end{lemma}
\proof
We prove the lemma by induction. Note that the lemma is true for $i=1$ as $\lambda_1=\lambda$ and $u$ is an eigenvector for $LR$ with eigenvalue $\lambda$.
Assume now that the lemma is true for each $j \: (1 \le j \le i-1)$. We will prove that it is true also for $i$. Note that $\alpha_{i+1} \ne 0$ and $Lu =0$, so 
by Corollary \ref{partII:cor0a} we have 
$$
  L R^i u = 
  \frac{c_{i+1} - \alpha_{i+1} }{\alpha_{i+1}} \Big( R^{i-1} L R u + \cdots + RLR^{i-1} u \Big) + \frac{t_{i+1} }{ \alpha_{i+1}} R^{i-1} u = 
$$
$$
\frac{c_{i+1} - \alpha_{i+1} }{ \alpha_{i+1}} \Big( \lambda_1 R^{i-1} u + \cdots + \lambda_{i-1} R^{i-1} u \Big) +\frac{t_{i+1} }{\alpha_{i+1}} R^{i-1} u = \lambda_{i} R^{i-1} u.
$$
The result follows. \qed

\begin{lemma}
\label{main:lem2}
With reference to Definition \ref{def:blank3} we have
$$
  \Big( c_{\ell}(\lambda_1 + \lambda_2 + \cdots + \lambda_{\ell-2}) + t_{\ell} \Big) R^{\ell-2} u = 0.
$$
\end{lemma}
\proof
As $\alpha_{\ell}=0$ and $Lu=0$, Corollary \ref{partII:cor0a} implies 
$$
  0 = c_{\ell} \Big( R^{\ell-2} L R u + \cdots + RLR^{\ell-2} u \Big) + t_{\ell} R^{\ell-2} u. 
$$
The result follows in view of Lemma \ref{main:lem1}. \qed

\begin{lemma}
\label{main:lem3}
With reference to Definition \ref{def:blank3} we have
$R^i u \ne 0$ for $0 \le i \le \ell-2$.
\end{lemma}
\proof
Assume to the contrary that $R^i u = 0$ for some $i \; (0 \le i \le \ell-2)$
and set $W = {\rm span} \{u, Ru, \ldots, R^{i-1} u\}$. We claim that $W$ is a $T$-module.
Observe that since $R \Es_j V \subseteq \Es_{j+1} V$ for $0 \le j \le D-1$ and $R\Es_DV=0$, $W$ is invariant 
under the action of $\Es_j$ for $0 \le j \le D$. Note also that $W$ is by definition invariant under   
the action of $R$. By Lemma \ref{main:lem1}, $W$ is also invariant under the action of $L$.
As $R+L = A$, it follows that $W$ is a $T$-module. Since $u \in U$, the endpoint of $W$ is 2.
As $W$ is a direct sum of irreducible $T$-modules, there exists an irreducible $T$-module $W'$
with endpoint $2$ which is contained in $W$. It follows that the diameter of $W'$ is less than or 
equal to $\ell-3$. But by \cite[Theorem 10.1]{Cu2}, the diameter of $W'$ is in $\{D-4,D-3,D-2\}$, and so 
$\ell-3 \ge D-4$, a contradiction. 
It follows that $R^i u \ne 0$ for $0 \le i \le \ell-2$. \qed

\begin{corollary}
\label{main:cor4}
With reference to Definition \ref{def:blank3} we have
$$
  c_{\ell}(\lambda_1 + \lambda_2 + \cdots + \lambda_{\ell-2}) + t_{\ell}=0.
$$
\end{corollary}
\proof
Immediate from Lemma \ref{main:lem2} and Lemma \ref{main:lem3}. \qed

\begin{theorem}
\label{thm:main}
With reference to Definition \ref{def:blank3} we have $|\Phi_2| \le 1$.  Furthermore, $\Delta_2 =0$.
\end{theorem}
\proof
By Lemma \ref{lem:loc1} and Lemma \ref{lem:loc2}, $\Phi_2$
consists of values $\eta$ for which $c_2 \eta + k$ is an eigenvalue of $LR\Es_2$ on $U$. Let $\eta \in \Phi_2$ such that $\lambda = c_2 \eta +k$. 
By Corollary \ref{main:cor4} and since $\lambda_i = \sigma_i \lambda + \rho_i$
for $1 \le i \le \ell-2$ we have
$$
  \lambda (\sigma_1 + \sigma_2 + \cdots + \sigma_{\ell-2}) = -\frac{t_{\ell} }{ c_{\ell}}-(\rho_1 + \rho_2 + \cdots + \rho_{\ell-2}). 
$$
As $\sigma_1 + \sigma_2 + \cdots + \sigma_{\ell-2} \ne 0$ by Lemma \ref{main:lem6}, we have
$$
  \lambda = -\frac{\frac{t_{\ell} }{ c_{\ell}}+(\rho_1 + \rho_2 + \cdots + \rho_{\ell-2}) }{ \sigma_1 + \sigma_2 + \cdots + \sigma_{\ell-2}}. 
$$
By definition the numbers $\rho_i, \sigma_i \; (1 \le i \le \ell-2)$ and $t_{\ell}$ depend only on the intersection numbers 
of $\G$ and on numbers $\alpha_i, \beta_i \; (3 \le i \le \ell-1)$. Therefore $\lambda$ (and hence also $\eta$) is uniquely determined, so $|\Phi_2| \le 1$. 
It follows from Lemma \ref{Del2} that $\Delta_2 =0$.  \qed

\section{Combinatorial property implies algebraic condition}

With reference to Notation \ref{blankII} assume $\Delta_2 > 0$. In this section we prove that
up to isomorphism $\G$ has exactly two irreducible $T$-modules 
with endpoint $2$, and they are both thin. Observe that by Theorem \ref{thm:main} we have $\alpha_i \ne 0$ for $3 \le i \le D-2$.

\begin{lemma}
\label{partII:lem0bb}
With reference to Notation \ref{blankII}, assume $\Delta_2 > 0$.  Let $W$ denote an irreducible $T$-module with 
endpoint $2$, and choose $u \in \Es_2 W$ which is an eigenvector for $LR$. 
Then $LR^{i-1} u \in {\rm Span}\{R^{i-2} u \}$ for $2 \le i \le D-2$.
\end{lemma}
\proof
We will prove the lemma by induction on $i$. Note that the statement of the lemma
is true for $i=2$ since $u$ is an eigenvector for $LR$. Assume $i \ge 3$.
The result now follows from \eqref{mat_dep1} using $\Es_2 u = u$, $L u = 0$, and the induction
hypothesis. \qed

\begin{corollary}
\label{partII:cor0bb}
With reference to Notation \ref{blankII}, assume $\Delta_2 > 0$.  Let $W$ denote an irreducible $T$-module with 
endpoint $2$, and choose $u \in \Es_2 W$ which is an eigenvector for $LR$. 
Then the following {\rm (i)--(iii)} hold.
\begin{itemize}
\item[{\rm (i)}]   $\Es_i A_{i-2} \Es_2 u = 1 / (c_{i-2} c_{i-3} \cdots c_1) R^{i-2} u$ $\quad (2 \le i \le D)$.
\item[{\rm (ii)}]  $\Es_i A_i \Es_2 u \in {\rm Span}\{R^{i-2} u\}$ $\quad (2 \le i \le D-2)$.
\item[{\rm (iii)}] $\Es_i A_{i+2} \Es_2 u \in {\rm Span}\{R^{i-2} u\}$ $\quad (2 \le i \le D-2)$.
\end{itemize}
\end{corollary}
\proof
(i) Immediate from Lemma \ref{LR:lem4}(i) using $\Es_2 u = u$.

\noindent
(ii) Immediate from Lemma \ref{LR:lem4}(ii) using $\Es_2 u = u$, $L u = 0$ and Lemma \ref{partII:lem0bb}.

\noindent
(iii) Note that by \cite[Corollary 9.3]{Cu2} we have 
$\Es_i A_{i+2} \Es_2 u = - \Es_i A_{i-2} \Es_2 u - \Es_i A_i \Es_2 u$. The result follows from
(i), (ii) above. \qed

\begin{corollary}
\label{partII:cor0bbb}
With reference to Notation \ref{blankII}, assume $\Delta_2 > 0$. Let $W$ denote an irreducible $T$-module with 
endpoint $2$, and choose $u \in \Es_2 W$ which is an eigenvector for $LR$. 
Then $LR^{D-2} u \in {\rm Span}\{R^{D-3} u \}$.
\end{corollary}
\proof
By Lemma \ref{LR:lem4}(i) and since $\Es_2 u = u$, we have $R^{D-2}u = c_1\cdots c_{D-2} \Es_D A_{D-2} \Es_2u$. 
By \cite[Theorem 9.4(iii)]{Cu2} we have
$L \Es_D A_{D-2} u = b_{D-1}\Es_{D-1} A_{D-3} u + R \Es_{D-2} A_D  u$.  Combining these facts, the result now follows from 
Corollary \ref{partII:cor0bb}(i),(iii). \qed

\begin{lemma}
\label{partII:lem0b}
With reference to Notation \ref{blankII}, assume $\Delta_2 > 0$.  Let $W$ denote an irreducible $T$-module with 
endpoint $2$ and diameter $d$. Choose $u \in \Es_2 W$ which is an eigenvector for $LR$. Then the
following is a basis for $W$:
\begin{equation}
\label{basis1}
  R^i u \qquad \qquad (0 \le i \le d).
\end{equation}
In particular, $W$ is thin.
\end{lemma}
\proof
We first show that $W$ is spanned by vectors \eqref{basis1}.
Let $W'$ denote the subspace of $V$ spanned by vectors \eqref{basis1} and note that $W' \subseteq W$.
We claim that $W'$ is $T$-invariant. Observe that since $R \Es_j V \subseteq \Es_{j+1} V$ for $0 \le j \le D-1$ and $R\Es_D V =0$, $W'$ is invariant 
under the action of $\Es_j$ for $0 \le j \le D$, and so $W'$ is $M^*$-invariant. 
By definition $W'$ is invariant under $R$. Note that by Lemma \ref{partII:lem0bb} and
Corollary \ref{partII:cor0bbb} $W'$ is invariant under $L$. Since $A=R+L$ and since $A$
generates $M$, $W'$ is $M$-invariant. The claim follows.
Hence $W'$ is a $T$-module and it is nonzero since $u \in W'$. By the irreducibility of $W$ we have $W'=W$.
Since $\Es_iW \ne 0$ for $2 \le i \le d+2$ we have $R^i u \ne 0$ for $0 \le i \le d$.
Observe also that since $R \Es_j V \subseteq \Es_{j+1} V$ for $0 \le j \le D-1$ and $R\Es_D V=0$, we have that $R^i u \; (0 \le i \le d)$ are linearly independent.
\qed

\begin{theorem}
\label{mainII}
With reference to Notation \ref{blankII} assume $\Delta_2 > 0$. Then $\G$ has up to isomorphism exactly two irreducible 
$T$-modules with endpoint $2$, and they are both thin.
\end{theorem}
\proof
Note that every irreducible $T$-module with endpoint $2$ is thin by Lemma \ref{partII:lem0b}.
We will now show that up to isomorphism $\G$ has exactly two irreducible 
$T$-modules with endpoint $2$.
By \cite[Theorem 4.2, Lemma 3.7, Theorem 3.8]{Maclean2}, the set $\Phi_2=\{\eta_{k+1}, \eta_{k+2}, \ldots, \eta_{k_2}\}$ has at most two elements. 

By \cite[Theorem 11.7]{Cu2}, the set $\Phi_2$ coincides with the set of local eigenvalues of the thin irreducible $T$-modules with endpoint 2.  
By Lemma \ref{isomequiv}, thin irreducible $T$-modules with endpoint $2$ 
are isomorphic if and only if they have the same local eigenvalue.  By Lemma \ref{Del2}, 
$\Phi_2$ has exactly two elements, and the result follows. \qed

\bigskip \noindent
We now present our main result.

\begin{theorem}
\label{mainresult}
Let $\Gamma$ denote a bipartite distance-regular graph with vertex set $X$, valency $k \ge 3$, and diameter $D \ge 4$.  
Then the following are equivalent.
\begin{itemize}
\item[{\rm (i)}] $\Delta_2>0$ and for $2 \le i \le D - 2$ there exist complex scalars 
$\alpha_i, \beta_i$ with the following property: for all $x, y, z \in X$ such that 
$\partial(x, y) = 2, \: \partial(x, z) = i, \: \partial(y, z) = i$ we have
$$
  \alpha_i + \beta_i  |\G(x) \cap \G(y) \cap \G_{i-1}(z)| = |\G_{i-1}(x) \cap \G_{i-1}(y) \cap \G(z)|.
$$
\item[{\rm (ii)}]
 For all $x \in X$ there exist up to isomorphism
exactly two irreducible modules for the Terwilliger algebra $T(x)$
with endpoint two, and these modules are thin.
\end{itemize}
\end{theorem}
\proof Immediate from Corollary \ref{partI:cor_mainI} and Theorem \ref{mainII}. \qed

\section{Scalars $\alpha_i$ and $\beta_i$ in terms of intersection numbers}

With reference to Notation \ref{blankII} assume that the equivalent conditions {\rm (i), (ii)}
of Theorem \ref{mainresult} hold. In this section we express scalars $\alpha_i$ and 
$\beta_i$ in terms of the intersection numbers of $\G$. To do this, we first need the following
definition.

\begin{definition}
\label{def:w}
With reference to Notation \ref{blank} and Definition \ref{def:D}, for all integers $i,j$ define a
vector $w_{ij}=w_{ij}(x,y)$ by
$$
  w_{ij} = \sum_{z \in D_j^i} {\hat z}.
$$
Observe that $w_{ij}=0$ unless $0 \le i,j \le D$ and $i+j$ is even, and that $\Vert w_{ij} \Vert^2 = p^2_{ij}$.

\smallskip \noindent
For all integers $i$ define vectors
$w_{ii}^+=w_{ii}^+(x,y)$ and $w_{ii}^-=w_{ii}^-(x,y)$ by
$$
  w_{ii}^+ = \sum_{z \in D_i^i} |\G_{i-1}(z) \cap D_1^1| \, {\hat z},
  \qquad \qquad
  w_{ii}^- = \sum_{z \in D_i^i} |\G(z) \cap D_{i-1}^{i-1}| \, {\hat z}.
$$
We observe $w_{ii}^+ = w_{ii}^- = 0$ unless $1 \le i \le D$. Furthermore,
$w_{11}^+ = w_{11}$, $w_{11}^- = 0$, $w_{DD}^+=c_2 w_{DD}$, $w_{DD}^-=kw_{DD}$, 
and $w_{22}^+ = w_{22}^-$.
\end{definition}
Theorem \ref{mainresult} can now be reformulated as follows.

\begin{theorem}
\label{mainresult-vector_form}
Let $\Gamma$ denote a bipartite distance-regular graph with vertex set $X$, valency $k \ge 3$, and diameter $D \ge 4$.  
Then the following are equivalent.
\begin{itemize}
\item[{\rm (i)}] $\Delta_2>0$ and for $2 \le i \le D - 2$ there exist complex scalars 
$\alpha_i, \beta_i$ with the following property: for all $x, y \in X$ such that 
$\partial(x, y) = 2$, we have
$$
  \alpha_i w_{ii} + \beta_i  w_{ii}^+ = w_{ii}^-,
$$
where vectors $w_{ii}, w_{ii}^+, w_{ii}^-$ are from Definition \ref{def:w}.
\item[{\rm (ii)}]
 For all $x \in X$ there exist up to isomorphism
exactly two irreducible modules for the Terwilliger algebra $T(x)$
with endpoint two, and these modules are thin.
\end{itemize}
\end{theorem}

\begin{lemma}
\label{bip:lem1}
{\rm (\cite[Lemma 4.1]{Mi2})}
With reference to Notation \ref{blank}, the following {\rm (i)--(iii)} hold.
\begin{itemize}
\item[{\rm (i)}]  $p^2_{i-2,i} = p^2_{i,i-2} = k_i c_i c_{i-1} / (k(k-1)) \;\;\; (2 \le i \le D)$;
\item[{\rm (ii)}] $p^2_{ii} = k_i (c_i(b_{i-1}-1)+b_i(c_{i+1}-1)) / (k(k-1))
                   \;\;\; (1 \le i \le D-1)$;
\item[{\rm (iii)}]  $p^2_{DD}=k_D (b_{D-1} - 1)/(k-1)$.
\end{itemize}
\end{lemma}

\begin{lemma}
\label{ort:lem1} 
{\rm (\cite[Lemma 7.1]{Mi2})}
With reference to Definition \ref{def:w} the following {\rm (i), (ii)} hold for 
$2 \le i \le D-1$.
\begin{itemize}
\item[{\rm (i)}]  $\la w_{ii}^+, w_{ii} \ra = k_i c_i (b_{i-1}-1) / k_2$;
\item[{\rm (ii)}] $\Vert w_{ii}^+ \Vert^2 = k_i c_i (c_2(b_{i-1}-1)-(c_2-1)b_i) / k_2$.
\end{itemize}
\end{lemma}

\begin{lemma}
\label{ort:lem2} 
With reference to Definition \ref{def:w} the following {\rm (i), (ii)} hold for 
$2 \le i \le D-1$.
\begin{itemize}
\item[{\rm (i)}]  $\la w_{ii}^-, w_{ii} \ra = c_i k_i (c_i b_{i-1} + c_{i-1} b_i - k) / (k(k-1))$;
\item[{\rm (ii)}] $\la w_{ii}^-, w_{ii}^+ \ra = k_i c_i (b_i(b_i-b_{i-1})+c_i(b_{i-1}-1)) / k_2$.
\end{itemize}
\end{lemma}
\proof
(i) We first claim that $\la c_i w_{ii} - w_{ii}^- , w_{ii} \ra = p^2_{i-1,i+1}(c_{i+1}-c_{i-1})$. 
To show this observe that 
$\la c_i w_{ii} - w_{ii}^-,w_{ii} \ra = \sum_{z \in D_i^i} |\G(z) \cap D_{i+1}^{i-1}|$.
Hence $\la c_i w_{ii} - w_{ii}^- , w_{ii} \ra$ is equal to the number of ordered pairs $(v,z)$,
where $v \in D_{i+1}^{i-1}$, $z \in D_i^i$, and $\partial(v,z)=1$. In order to find this number, we 
fix $v \in D_{i+1}^{i-1}$ and observe $|\G(v) \cap D_i^i|=c_{i+1} - c_{i-1}$.
The claim follows since there are exactly $p^2_{i-1, i+1}$ vertices in $D_{i+1}^{i-1}$.
The result now follows from the above claim, Lemma \ref{bip:lem1}, and since $\Vert w_{ii} \Vert^2 = p^2_{ii}$.

\medskip \noindent
(ii) We first claim that $\la (b_i-c_i) w_{ii} + w_{ii}^-, w_{ii}^+ \ra = c_2 c_i (p^1_{i+1,i}-p^2_{i+1,i-1})$.
To show this observe that
$\la (b_i - c_i)w_{ii} + w_{ii}^-,w_{ii}^+ \ra = 
 \sum_{z \in D_i^i} |\G(z) \cap D_{i+1}^{i+1}| |\G_{i-1}(z) \cap D_1^1|$. 
Hence $\la (b_i - c_i)w_{ii} + w_{ii}^-, w_{ii}^+ \ra$ is equal to the number of ordered triples 
$(v,z,u)$, where $v \in D_1^1$, $z \in D_i^i$, $u \in D_{i+1}^{i+1}$, $\partial(v,z) = i-1$ and 
$\partial(z,u)=1$. Therefore, $\partial(v,u) = i$. 
Now pick $v \in D_1^1$ and note that we have $c_2$ choices for this. 
Observe that $\G_{i+1}(x) \cap \G_i(v) \subseteq D^{i+1}_{i-1} \cup D_{i+1}^{i+1}$.
Moreover, every vertex in $D_{i-1}^{i+1}$ is at distance $i$ from $v$. Therefore, since $u$ is in $D_{i+1}^{i+1}$ and $\partial(v,u) = i$,
we have $p^1_{i+1,i} - p^2_{i+1,i-1}$ choices for $u$.
Finally, for every such pair $(v,u)$ there are exactly $c_i$ vertices $z$ which are at distance 
$i-1$ from $v$ and at distance 1 from $u$. Observe that all these vertices must be in $D_i^i$. This proves the claim.
The result now follows from the claim, Lemma \ref{ort:lem1}(i), Lemma \ref{bip:lem1} and since $p^1_{i+1,i} = b_i k_i / k$. \qed

\begin{theorem}
\label{thm:scalars}
  Let $\Gamma$ denote a bipartite distance-regular graph with vertex set $X$, valency $k \ge 3$, and diameter $D \ge 4$.  
Assume that the equivalent conditions {\rm (i), (ii)} of Theorem \ref{mainresult-vector_form} hold.
Then for $2 \le i \le D-2$, $\Delta_i > 0$ and the following hold:
$$
\alpha_i = \frac{c_i (c_i-1) (b_{i-1}- c_2) - c_i c_{i-1} (b_i-1) (c_2-1) }{c_2 \Delta_i}
$$
and 
$$
\beta_i = \frac{c_i (c_{i+1}-c_i) (b_{i-1}-1) - b_i (c_{i+1}-1) (c_i-c_{i-1}) }{c_2 \Delta_i},
$$
where scalars $\Delta_i$ are from Definition \ref{def:delta}.
\end{theorem}
\proof
By Lemma \ref{Del2}, $\Delta_i \ge 0$.
Suppose $\Delta_i =0$.  Then by \cite[Theorem 13]{Cu0}, for all $x,y,z\in X$ with $\partial(x,y)=2$, $\partial(x,z)=i$, $\partial(y,z)=i$, the number $|\Gamma(x) \cap \Gamma(y) \cap \Gamma_{i-1}(x)|$ is independent of $x,y,z$.  By Theorem \ref{partI:thm_mainI}, we find (\ref{eq:recurrence}) is satisfied with $\gamma_i =0$, which contradicts 
Theorem \ref{thm:gamma}.  Hence $\Delta_i >0$.
To obtain the formulae for $\alpha_i$ and $\beta_i$, take the inner product of 
the vector equation in Theorem \ref{mainresult-vector_form} with $w_{ii}$ and $w_{ii}^+$,
and then solve thus obtained system of linear equations using Lemma \ref{ort:lem1} and Lemma \ref{ort:lem2}. \qed


\begin{thebibliography}{99}

{\small

\bibitem{BI} E.\ Bannai, and T.\ Ito, {\em Algebraic Combinatorics I: Association schemes},
Benjamin-Cummings Lecture Note 58, Menlo Park, 1984.

\bibitem{BCN} A.\ E.\ Brouwer, A.\ M.\ Cohen, and A.\ Neumaier,
 {\em Distance-Regular Graphs}, Springer-Verlag, Berlin, 1989.
   
\bibitem{Cu0} B.\ Curtin,
2-homogeneous bipartite distance-regular graphs, {\it Discrete Math.} 
{\bf 187} (1998), 39-70.

\bibitem{Cu1} B.\ Curtin, Bipartite distance-regular graphs, Part I,
 {\it Graphs Combin.} {\bf 15} (1999), 143-158.
 
\bibitem{Cu2} B.\ Curtin, Bipartite distance-regular graphs, Part II,
 {\it Graphs Combin.} {\bf 15} (1999), 377-391.
 
\bibitem{Cu3} B.\ Curtin,  The local structure of a bipartite distance-regular graph,
 {\it European J. Combin.} {\bf 20} (1999), 739-758.
 
\bibitem{Cu4} B.\ Curtin,  Almost $2$-homogeneous bipartite distance-regular graphs,
 {\it European J. Combin.} {\bf 21} (2000), 865-876.

\bibitem{egge1} E.\ Egge, A generalization of the Terwilliger algebra,
{\it J.\ Algebra} {\bf 233} (2000), 213-252.

\bibitem{Go1} J.\ Go, The Terwilliger algebra of the Hypercube,
{\it European J. Combin.} {\bf 23} (2002), 399-429.

 
\bibitem{MT} M.~S.~MacLean and P.~Terwilliger, Taut distance-regular graphs and the 
subconstituent algebra, {\it Discrete Math.} {\bf 306} (2006), 1694-1721.

\bibitem{MT2}
M.~S.~MacLean and P.~Terwilliger, The subconstituent 
algebra of a bipartite distance-regular graph:  thin modules with endpoint two, {\it Discrete Math.} {\bf 308} (2008), 1230-1259.


\bibitem{Maclean2} M.~S.~MacLean, The local eigenvalues of a bipartite distance-regular graph, 
{\it European J. Combin.} {\bf 45} (2015), 115-123.

\bibitem{Mi2}  \v{S}.\ Miklavi\v{c}, On bipartite $Q$-polynomial distance-regular graphs,
{\it European J. Combin.} {\bf 28} (2007), 94-110.

\bibitem{No1} K.\ Nomura, Homogeneous graphs and regular near polygons,
 {\it J. Combin. Theory Ser. B} {\bf 60} (1994), 63-71.
 

\bibitem{ter1} P.\ Terwilliger, The subconstituent algebra of an association scheme, Part I,
 {\em J.\ Algebr.\ Combin.} {\bf 1} (1992), 363-388.
 

}

\end{thebibliography}
\end{document}